\documentclass[10pt]{amsart}
\usepackage{latexsym,amssymb,amsmath,verbatim}

\theoremstyle{plain}
  \newtheorem{theorem}{Theorem}
  \newtheorem{proposition}[theorem]{Proposition}
  \newtheorem{lemma}[theorem]{Lemma}
  \newtheorem{corollary}[theorem]{Corollary}
  \newtheorem{conjecture}[theorem]{Conjecture}
\theoremstyle{definition}
  \newtheorem{definition}[theorem]{Definition}

  \newtheorem{question}[theorem]{Question}
 \theoremstyle{remark}

\numberwithin{equation}{section}

\newcommand\qbin[3]{\left[\begin{matrix} #1 \\ #2 \end{matrix} \right]_{#3}}

\newcommand\inv{\operatorname{inv}}
\newcommand\sgn{\operatorname{sgn}}
\newcommand\triv{\mathbf{1}}
\newcommand\Tr{\operatorname{Tr}}
\newcommand\rank{\operatorname{rank}}
\newcommand\wt{\operatorname{wt}}

\newcommand\induce{\uparrow}
\newcommand\restrict{\downarrow}
\newcommand\inflate{\Uparrow}
\newcommand\semidirect{\ltimes}

\newcommand\St{\operatorname{St}}
\newcommand\ch{\operatorname{ch}}

\newcommand\FF{\mathbb{F}}

\newcommand\GG{\mathbb{G}}
\newcommand\ZZ{\mathbb{Z}}
\newcommand\NN{\mathbb{N}}
\newcommand\RR{\mathbb{R}}
\newcommand\CC{\mathbb{C}}

\newcommand\UUU{\mathcal{U}}

\newcommand\symm{\mathfrak{S}}

\begin{document}
\title[The negative $q$-binomial]{The negative $q$-binomial}
\thanks{Version of August 8, 2011}

\author{Shishuo Fu}
\email{fu@math.psu.edu}
\address{Department of Mathematics\\
The Pennsylvania State University\\
University Park, PA 16801}

\author{Victor Reiner}
\email{reiner@math.umn.edu}
\address{School of Mathematics\\
University of  Minnesota\\
Minneapolis, MN 55455}

\author{Dennis Stanton}
\email{stanton@math.umn.edu}
\address{School of Mathematics\\
University of  Minnesota\\
Minneapolis, MN 55455}

\author{Nathaniel Thiem}
\email{thiemn@colorado.edu}
\address{Department of Mathematics\\
University of Colorado at Boulder\\
Boulder, CO 80309}

\begin{abstract} Interpretations for the
$q$-binomial coefficient evaluated at $-q$ are
discussed.  A $(q,t)$-version is established,
including an instance of a cyclic sieving phenomenon
involving unitary spaces.
\end{abstract}

\maketitle

\section{The $q$-binomial}

 The $q$-binomial coefficient is defined for
integers $k$ and $n$, with $0\le k\le n$, and an indeterminate $q$ by
\begin{equation}
\label{q-binomial-defn}
\qbin{n}{k}{q}=
\frac{(q)_n}{(q)_k(q)_{n-k}}
\end{equation}
where $(q)_n=(1-q^1)(1-q^2)\cdots(1-q^n)$.
It is well-known \cite[p. 39]{A} that the 
$q$-binomial coefficient is a polynomial in
$q$ with non-negative integer coefficients
\begin{equation}
\label{gf}
\qbin{n}{k}{q}=\sum_{\omega\in \Omega_{n,k}} q^{\inv(\omega)},
\end{equation}
where 
$\Omega_{n,k}$ is the set of words 
$\omega=(\omega_1,\ldots,\omega_n)$ in $\{0,1\}^n$
having $k$ ones and $n-k$ zeroes, 
and $\inv(\omega)$ is the number of {\it inversions} in $\omega$, 
that is, pairs $(i,j)$ with $1 \leq i < j \leq n$ and
$\omega_i=1, \omega_j=0$; see \cite[p. 40]{A}.
When $q$ is a prime power, the $q$-binomial coefficient 
\eqref{q-binomial-defn} is an integer
counting the number of  $k$-dimensional spaces in the $n$-dimensional vector space
$\FF_q^n$ over the field $\FF_q$. 

Section ~\ref{definition-section} combinatorially
interprets the $q$-binomial coefficient when $q$ is a 
negative integer (Theorem~\ref{mainth}),
while Section~\ref{q-t-section} establishes a 
positivity theorem for a $(q,t)$-analogue when $q$ is negative 
(Theorem ~\ref{lastthm}).  Section~\ref{Ennola-section}
provides a different interpretation for the negative $q$-binomial,
counting unitary subspaces, and related to {\it Ennola duality} for
finite unitary groups.  Section~\ref{CSP-section} 
proves a cyclic sieving phenomenon
involving the $(q,t)$-analogue at negative $q$ and unitary subspaces. 
Section~\ref{remarks-section} collects some
remarks and remaining questions suggested by these results.

\section{The negative $q$-binomial}
\label{definition-section}

Let $q\ge2$ be an integer and define
\begin{equation}
\label{negative-q-binomial-definition}
\qbin{n}{k}{q}':=(-1)^{k(n-k)}\qbin{n}{k}{-q}.
\end{equation}
\noindent
It is not hard to derive from the
the product expression in \eqref{q-binomial-defn} that this 
primed $q$-binomial will be positive, and it follows
from \eqref{gf} that it is an integer, with
\begin{equation}
\label{ineq}
\qbin{n}{k}{q}'\le \qbin{n}{k}{q}.
\end{equation}

\noindent
Our main result Theorem~\ref{mainth} 
is an analogue of \eqref{gf} for the primed
$q$-binomial coefficient that clearly demonstrates \eqref{ineq}. 
It expands the primed
$q$-binomial coefficient as a sum over
words $\omega$ in a 
subset $\Omega'_{n,k} \subset \Omega_{n,k}$, with weights 
$\wt(\omega)$ satisfying 
$1\le \wt(\omega)\le q^{\inv(\omega)}$ for $q \geq 2$.
This subset $\Omega'_{n,k}$ and
weight $\wt(\omega)$ come from a
{\it pairing} algorithm explained next.

\vskip.1in
\noindent
{\bf Definition.}
Given $\omega=(\omega_1,\ldots,\omega_n)$ 
in $\Omega_{n,k}$, pair some of its adjacent entries 
$(\omega_{i},\omega_{i+1})$, and leave others unpaired,
according to the following recursive rule:
\begin{enumerate}
\item[$\bullet$]
When $n=1$, leave the unique letter $\omega_1$ in $\omega$ unpaired.
\item[$\bullet$]
When $n \geq 2$ and $k$ is odd, pair 
the first two entries $\omega_1,\omega_2$,
and recursively pair the remaining word $(\omega_3,\ldots,\omega_n)$.
\item[$\bullet$]
When $n \geq 2$ and $k$ is even, 
leave the first entry $\omega_1$ unpaired, 
and recursively pair the remaining word $(\omega_2,\omega_3,\ldots,\omega_n)$.
\end{enumerate}
Two examples of  pairings for words, with pairings
indicated by underlining, are
$$
\begin{aligned}
\omega^{(1)}=&\ \underline{01}\ \underline{1}\ \underline{00}\ \underline{10}
\ \underline{1}\ \underline{01}\\
\omega^{(2)}=&\  \underline{1}\ \underline{10}\ \underline{0}\
\underline{0}\ \underline{1}\ \underline{00} \ \underline{1}.
\end{aligned}
$$
\noindent
Define 
$$
\Omega'_{n,k}:=\{ \omega \in \Omega_{n,k}: 
                 \omega \text{ has no paired }\underline{01}\}.
$$
For example, $\omega^{(2)}$ lies in $\Omega'_{n,k}$, but
$\omega^{(1)}$ does not.
For $\omega\in\Omega'_{n,k}$, define
$$
\begin{array}{rl}
p(\omega) &:=\text{number of }\underline{10}\text{ pairs in }\omega  \\
 &   \\
a(\omega)&:=\inv(\omega)-p(\omega) \\
 &   \\
\wt(\omega)&:=q^{a(\omega)} (q-1)^{p(\omega)}.
\end{array}
$$
\noindent
Note that $a(\omega) \geq 0$ since each pair $\underline{10}$ contributes
at least $1$ to the value of $\inv(\omega)$.  In fact, it is helpful
to think of $a(\omega)$ as a 
perturbation of the inversion statistic
$\inv(\omega)$ as follows:  each occurrence of $1$ in $\omega$
would normally contribute to $\inv(\omega)$ the number of zeroes to its right,
but when this $1$ occurs in a pair $\underline{10}$ it contributes
one fewer than usual to $a(\omega)$. 

Note also that
$1 \leq \wt(\omega) \leq q^{\inv(\omega)}$ for $q \geq 2.$
This brings us to the main result. 

\begin{theorem} For $0\le k\le n$, one has 
\label{mainth}
$$
\qbin{n}{k}{q}'=\sum_{\omega\in\Omega'_{n,k}} \wt(\omega)
=\sum_{\omega\in\Omega'_{n,k}} q^{a(\omega)} (q-1)^{p(\omega)}.
$$
\end{theorem}

For example, if $(n,k)=(5,2)$, one has this table of data and calculation:
$$
\begin{tabular}{|c|c|}\hline
$\omega \in \Omega'_{5,2}$ & $\wt(\omega)$ \\ \hline\hline
$\underline{0}\ \underline{0}\ \underline{0}\ \underline{1}\ \underline{1}$
&$1$\\\hline
$\underline{0}\ \underline{0}\ \underline{1}\ \underline{10}$
&$q(q-1)$\\\hline
$\underline{0}\ \underline{1}\ \underline{00}\ \underline{1}$
&$q^2$\\\hline
$\underline{0}\ \underline{1}\ \underline{10}\ \underline{0}$
&$q^3(q-1)$\\\hline 
$\underline{1}\ \underline{00}\ \underline{10}$
&$q^3(q-1)$\\\hline 
$\underline{1}\ \underline{10}\ \underline{0}\ \underline{0}$
&$q^5(q-1)$\\\hline
\end{tabular}
\qquad 
\vspace{.1in}
\begin{aligned}
&1+q(q-1)+q^2+q^3(q-1) +q^3(q-1)+q^5(q-1)\\
&=q^6-q^5+2q^4-2q^3+2q^2-q+1\\
&=(-1)^{2 \cdot (5-2)}\qbin{5}{2}{-q}= \qbin{5}{2}{q}'.
\end{aligned}
$$
\begin{proof}[Proof of Theorem \ref{mainth}]
Induct on $n$, with easily verified base cases $n=0,1$.
In the inductive step, we use this
$q$-Pascal recurrence \eqref{q-Pascal}, and 
its iterate \eqref{iterated-q-Pascal}
\begin{equation}
\label{q-Pascal}
\qbin{n}{k}{q} =
\qbin{n-1}{k}{q} +
q^{n-k} \qbin{n-1}{k-1}{q},
\end{equation}
\begin{equation}
\label{iterated-q-Pascal}
\qbin{n}{k}{q}=
\qbin{n-2}{k}{q} +
q^{n-k-1} (q+1) \qbin{n-2}{k-1}{q} +
q^{2(n-k)} \qbin{n-2}{k-2}{q}.
\end{equation}
When $k$ is even, replacing $q \mapsto -q$ in \eqref{q-Pascal} gives
\begin{equation}
\label{negative-q-Pascal}
\qbin{n}{k}{q}'=\qbin{n-1}{k}{q}' + q^{n-k} \qbin{n-1}{k-1}{q}'.
\end{equation}
Since $k$ is even, the leading entry $\omega_1$ will be either be
an unpaired $\underline{0}$ or $\underline{1}$.  In either case,
$\omega_1$ contributes $0$ to $p(\omega)$.  If $\omega_1=\underline{0}$
it contributes $0$ to $a(\omega)$, and corresponds to the first
summand on the right of \eqref{negative-q-Pascal}, while if
$\omega_1=\underline{1}$ it contributes $n-k$ to $a(\omega)$, and 
corresponds to the second summand on the right of \eqref{negative-q-Pascal}.

When $k$ is odd, replacing $q \mapsto -q$ in \eqref{iterated-q-Pascal} gives
\begin{equation}
\label{negative-iterated-q-Pascal}
\qbin{n}{k}{q}'=
q^0(q-1)^0 \qbin{n-2}{k}{q}'+
q^{n-k-1} (q-1)^1 \qbin{n-2}{k-1}{q}' +
q^{2(n-k)}(q-1)^0 \qbin{n-2}{k-2}{q}'.
\end{equation}
Since $k$ is odd, $(\omega_1,\omega_2)$ will be paired, and
since $\omega$ lies in $\Omega'_{n,k}$, the pair $(\omega_1,\omega_2)$
takes one of
the three forms $\underline{00}$, $\underline{10}$, or $\underline{11}$,
\begin{enumerate}
\item[$\bullet$]
contributing $0$, $n-k-1$, or $2(n-k)$, respectively, to $a(\omega)$, 
\item[$\bullet$]
contributing $0$, $1$, or $0$, respectively, to $p(\omega)$, and
\item[$\bullet$]
leaving $k$ ones, $k-1$ ones, or $k-2$ ones, respectively, in $(\omega_3,\ldots,\omega_n)$.
\end{enumerate}
Thus the three forms correspond to the three summands of \eqref{negative-iterated-q-Pascal}
\end{proof}

We note that one can reformulate Theorem \ref{mainth}
as an expansion of the $q$-binomial coefficient, with
no negative signs, as follows.
\begin{corollary} If $0\le k\le n$,
\label{maincor}
$$
\qbin{n}{k}{q}=
\sum_{\omega\in\Omega'_{n,k}}
q^{a(\omega)}
(q+1)^{p(\omega)}.
$$
\end{corollary}
\begin{proof} Setting $q \mapsto -q$ in Theorem~\ref{mainth} 
and multiplying by $(-1)^{k(n-k)}$ gives
$$
\qbin{n}{k}{q}=
\sum_{\omega\in\Omega'_{n,k}}
(-1)^{a(\omega)+p(\omega)+k(n-k)} q^{a(\omega)} (q+1)^{p(\omega)}.
$$
Comparing with the corollary,
it suffices to show that for each $\omega$ in 
$\Omega'_{n,k}$ one has the following parity condition:
\begin{equation}
\label{parity-condition}
\inv(\omega) = a(\omega)+p(\omega) \quad \equiv \quad k(n-k)\pmod 2.
\end{equation}
This can be checked via induction on $n$ 
using the recursive definition of  $\Omega'_{n,k}$:
\vskip.1in
\noindent
{\sf Case 1.} $k$ is even.
\newline  
If $\omega=\underline{0} \omega'$, then 
$
\inv(\omega)=\inv(\omega')
\equiv 
k(n-1-k)\equiv k(n-k) \bmod 2.
$
\newline
If $\omega=\underline{1}\omega'$, then
$
\inv(\omega)=n-k+\inv(\omega')
\equiv 
n-k+(k-1)(n-k)\equiv k(n-k) \bmod 2.
$
\vskip.1in
\noindent
{\sf Case 2.} $k$ is odd.  
\newline
If $\omega=\underline{00}\omega'$ or $\omega=\underline{11}\omega'$, then
$
\inv(\omega) \equiv \inv(\omega') 
\equiv 
k(n-k-2) \equiv k(n-k) \bmod 2.
$
\newline
If $\omega=\underline{10}\omega'$, then
$$
\inv(\omega)=n-k+\inv(\omega')
\equiv
n-k+(k-1)(n-1-k)
\equiv k(n-k) \bmod 2.
$$
\end{proof}


%

\section{The $(q,t)$-binomial at negative $q$}
\label{q-t-section}

In \cite[p. 43]{RSW} the authors consider a certain
$(q,t)$-analogue of the $q$-binomial: a polynomial in $t$ 
with positive integer coefficients, depending upon
a positive integer $q$, and whose limit as $t$ goes to $1$
is the $q$-binomial \cite[Corollary 3.2]{RS}. 
Here we first review the definition and interpretation
of this $(q,t)$-binomial, and then 
establishing a positivity result for it (Theorem \ref{lastthm})
when $q$ is a {\it negative} integer.

This $(q,t)$-version of the binomial coefficient
is defined by (see  \cite{RS} or \cite[p. 43]{RSW})
\begin{equation}
\label{q-t-binomial-definition}
\qbin{n}{k}{q,t}:=\prod_{i=1}^k \frac{1-t^{q^n-q^{i-1}}}
{1-t^{q^k-q^{i-1}}}.
\end{equation}
When $q$ a positive integer, one can show\footnote{It follows, 
e.g., by iterating the $(q,t)$-Pascal recurrence 
\eqref{qt-Pascal} from the proof of Theorem~\ref{lastthm}.}
that this turns out to be a polynomial in $t$ with 
nonnegative integer coefficients.
One can easily check from the above definition that 
its degree in $t$ is $k(q^n-q^k)$, and that
\begin{equation}
\label{coefficient-symmetry}
t^{k(q^n-q^k)} \qbin{n}{k}{q,t^{-1}}
=
 \qbin{n}{k}{q,t}
\end{equation}
so that its coefficient sequence will be symmetric about the power
$t^{\frac{1}{2} k(q^n-q^k)}$.

When $q$ is a prime power, and hence the order of a finite
field $\FF_q$, these properties of the
$(q,t)$-binomial follow from its interpretation
as the Hilbert series for a certain graded ring, 
briefly reviewed here; see \cite[\S 9]{RSW} or \cite{RS}
for more detail. One starts with a 
$S:=\FF_q[x_,\ldots,x_n]$ a polynomial algebra
on which $G:=GL_n(\FF_q)$ acts by linear substitutions of variables.
One has nested subalgebras $S^G \subset S^P$ of 
$G$-invariant polynomials $S^G$, and $P$-invariant polynomials $S^P$, where
$P$ is the parabolic subgroup of $G$ that stabilizes a particular $k$-dimensional
subspace of $\FF_q^n$.  Then the $(q,t)$-binomial coefficient in \eqref{q-t-binomial-definition}
is the Hilbert series in the variable $t$ for the quotient ring $S^P/(S^G_+)$, in which
$(S^G_+)$ denotes the ideal of  $S^P$ generated by the $G$-invariant polynomials
of strictly positive degree.

\begin{theorem}
\label{lastthm}
When $q\le -2$ is a negative integer, the $(q,t)$-binomial
defined as a rational function in \eqref{q-t-binomial-definition}
lies in $(-1)^{k(n-k)} \NN[t,t^{-1}]$, meaning that
it is a Laurent polynomial in $t$ whose nonzero
coefficients all have sign $(-1)^{k(n-k)}$.

Furthermore, its coefficient sequence is
symmetric about $t^{\frac{1}{2} k(q^n-q^k)}$,
with monic coefficients on its smallest and largest
powers of $t$, which are the following powers:
$$
\begin{array}{rcccll}
\{& t^0 &,& t^{k(q^n-q^k)} &\} & \text{ if }n,k\text{ are both even}, \\
\{& t^{k(q^n-q^k)} &,& t^0 &\} & \text{ if }n,k\text{ are both odd}, \\
\{& t^{kq^n-\frac{1-q^k}{1-q}} &,& t^{-kq^k+\frac{1-q^k}{1-q}} &\} 
   & \text{ if }n\text{ is odd and }k\text{ is even}, \\
\{& t^{-kq^k+\frac{1-q^k}{1-q}} &,& t^{kq^n-\frac{1-q^k}{1-q}} &\}
   & \text{ if }n\text{ is even and }k\text{ is odd}. \\
\end{array}
$$
\end{theorem}

\begin{proof} 
The main issue is proving that this $(q,t)$-binomial
with $q \le -2$ is a Laurent polynomial of the appropriate
sign.  Given this, the coefficient symmetry follows from the 
validity of \eqref{coefficient-symmetry} 
for any integer value of $q$.  The last assertion
of the theorem follows by symmetry, after examining in each case 
the beginning of the Laurent expansion of the product on the
right side of \eqref{q-t-binomial-definition}.  We omit
the details.

The first assertion is proven by induction on $n$,
as in the proof of Theorem \ref{mainth}.
The base cases $n=0,1$ are again easily verified.  
In the inductive step, one proceeds in two cases, 
based on the parity of $n-k$.

\vskip.1in
\noindent
{\sf Case 1.} $n-k$ is even.
\newline
We use the analogue \cite[Prop. 4.1]{RS} of \eqref{q-Pascal} 
with $k$ replaced by $n-k$: 
\begin{equation}
\label{qt-Pascal}
\qbin{n}{k}{q,t} = \qbin{n-1}{k-1}{q,t^q}+t^{q^{k}-1}\prod_{i=0}^{k-1}\dfrac{1-t^{q^{k+1}-q^{i+1}}}{1-t^{q^{k}-q^{i}}}\qbin{n-1}{k}{q,t^q}
\end{equation}

\noindent
We check both summands on the right of \eqref{qt-Pascal} lie
in $(-1)^{k(n-k)}\NN[t,t^{-1}] = \NN[t,t^{-1}]$.
By induction on $n$, the first term 
lies in $(-1)^{(k-1)(n-k)} \NN[t,t^{-1}] = \NN[t,t^{-1}]$.

For the second term, again by induction on $n$,
its $(q,t)$-binomial factor lies 
in $(-1)^{k(n-k-1)}\NN[t,t^{-1}] = (-1)^k \NN[t,t^{-1}]$.
It then suffices to verify that
the product over $i=0,1,\ldots,k-1$
in the second term has each of its $k$ factors
lying in $(-1)^1 \NN[t,t^{-1}]$.  To verify this, let $m:=q^k-q^i$
and then this factor can be rewritten
\begin{equation}
\label{negative-geometric-series}
\begin{aligned}
\dfrac{1-t^{q^{k+1}-q^{i+1}}}
      {1-t^{q^{k}-q^{i}}}
&=\dfrac{1-t^{qm}}
       {1-t^m}
=-t^{qm}\dfrac{1-t^{-qm}}
       {1-t^m}\\
&=-t^{qm} \left(1+t^m +t^{2m}+\cdots+t^{(-q-1)m} \right).
\end{aligned}
\end{equation}
This lies in $(-1)^1 \NN[t,t^{-1}]$.  Hence
the second term of \eqref{qt-Pascal} lies in $\NN[t,t^{-1}]$.

\vskip.1in
\noindent
{\sf Case 1.} $n-k$ is odd.
\newline
We use the analogue of \eqref{iterated-q-Pascal} with $k$ replaced by $n-k$:
\begin{equation}
\label{iterated-qt-Pascal}
\qbin{n}{k}{q,t} = A+B+C+D 
\end{equation}
$$
\begin{aligned}
\text{ where }\quad 
A:=& \qbin{n-2}{k-2}{q,t^{q^2}} \\
B:=& t^{q^{k}-q}\prod_{i=0}^{k-2}\dfrac{1-t^{q^{k+1}-q^{i+2}}}{1-t^{q^{k}-q^{i+1}}}\qbin{n-2}{k-1}{q,t^{q^2}}\\
C:=& t^{q^{k}-1}\prod_{i=0}^{k-1}\dfrac{1-t^{q^{k+1}-q^{i+1}}}{1-t^{q^{k}-q^{i}}}\qbin{n-2}{k-1}{q,t^{q^2}}\\
D:=& t^{q^{k+1}+q^{k}-q-1}\prod_{i=0}^{k-1}\dfrac{1-t^{q^{k+2}-q^{i+2}}}{1-t^{q^{k}-q^{i}}}\qbin{n-2}{k}{q,t^{q^2}}
\end{aligned}
$$
\noindent
The last term $D$ is easy to deal with alone.  By induction
on $n$, its $(q,t)$-binomial factor lies in 
$(-1)^{k(n-2-k)} \NN[t,t^{-1}] =
(-1)^{k(n-k)} \NN[t,t^{-1}]$.  We claim that 
each factor for $i=0,1,\ldots,k-1$ within the product inside
$D$ lies in $\NN[t,t^{-1}]$, since it can be expressed
$$
\dfrac{1-t^{q^{k+2}-q^{i+2}}}{1-t^{q^{k}-q^{i}}}=
\dfrac{1-t^{q^2m}}{1-t^{m}}=
1+t^m + t^{2m} +\cdots + t^{(q^2-1)m}.
$$
where $m:=q^k-q^i$ as before.

The factors in $B$ and $C$ which correspond to 
$q^{n-k-1}(q-1)^1$ in  \eqref{negative-iterated-q-Pascal}, 
are no longer Laurent polynomials in $t$ with 
non-negative coefficients. Thus more care 
must be taken to prove that $A+B+C$ lies in
$(-1)^{k(n-k)} \NN[t,t^{-1}]$ for $q \leq -2$.  
We start by combining common factors in $B$ and $C$:
$$
B+C
=\prod_{i=1}^{k-1}\dfrac{1-t^{q^{k+1}-q^{i+1}}}{1-t^{q^{k}-q^{i}}}
    \qbin{n-2}{k-1}{q,t^{q^2}}
     \Big(t^{q^{k}-q}+t^{q^{k}-1} \dfrac{1-t^{q^{k+1}-q}}{1-t^{q^{k}-1}}\Big) 
$$
and rewrite this parenthesized factor within $B+C$ as follows:
$$
\begin{aligned}
&t^{q^{k}-q}+t^{q^{k}-1} \dfrac{1-t^{q(q^{k}-1)}}{1-t^{q^{k}-1}} \\
&=t^{q^{k}-q}-t^{q^{k}-1}\big(t^{(-1)(q^{k}-1)}+t^{(-2)(q^{k}-1)}+\cdots+t^{q(q^{k}-1)}\big)\\
&=(t^{q^{k}-q}-1)- h\\
\text{ where } h
  &:=\big(t^{(-1)(q^{k}-1)}+t^{(-2)(q^{k}-1)}+\cdots+t^{(q+1)(q^{k}-1)}\big).\\
\end{aligned}
$$
One can also rewrite the other two factors appearing in $B+C$:
$$
\begin{aligned}
\prod_{i=1}^{k-1}\dfrac{1-t^{q^{k+1}-q^{i+1}}}{1-t^{q^{k}-q^{i}}} \qbin{n-2}{k-1}{q,t^{q^2}} 
&=\prod_{i=1}^{k-1}\dfrac{1-t^{q^{k+1}-q^{i+1}}}{1-t^{q^{k}-q^{i}}} 
\prod_{i=0}^{k-2}\dfrac{1-t^{q^{n}-q^{i+2}}}{1-t^{q^{k+1}-q^{i+2}}} \\
&=\prod_{i=0}^{k-2}\dfrac{1-t^{q^{n}-q^{i+2}}}{1-t^{q^{k}-q^{i+1}}} \quad \text{( via telescoping )}\\
&=\dfrac{1-t^{q^{n}-q^{k}}}{1-t^{q^{k}-q}}
    \prod_{i=0}^{k-3}\dfrac{1-t^{q^{n}-q^{i+2}}}{1-t^{q^{k}-q^{i+2}}}\\
&=\dfrac{1-t^{q^{n}-q^{k}}}{1-t^{q^{k}-q}}
    \qbin{n-2}{k-2}{q,t^{q^2}}.
\end{aligned}
$$
Therefore $A+B+C$ equals
$$
\begin{aligned}
&\qbin{n-2}{k-2}{q,t^{q^2}} + 
         \prod_{i=1}^{k-1}\dfrac{1-t^{q^{k+1}-q^{i+1}}}{1-t^{q^{k}-q^{i}}}
         \qbin{n-2}{k-1}{q,t^{q^2}} \cdot ((t^{q^{k}-q}-1)- h)  \\
&=\qbin{n-2}{k-2}{q,t^{q^2}} + 
         \dfrac{1-t^{q^{n}-q^{k}}}{1-t^{q^{k}-q}} 
           \qbin{n-2}{k-2}{q,t^{q^2}}(t^{q^{k}-q}-1)\\
&\qquad \qquad 
-       \prod_{i=1}^{k-1}\dfrac{1-t^{q^{k+1}-q^{i+1}}}{1-t^{q^{k}-q^{i}}} 
         \qbin{n-2}{k-1}{q,t^{q^2}} \cdot h  \\
&=  \qbin{n-2}{k-2}{q,t^{q^2}} 
      \left( 1 + \dfrac{1-t^{q^{n}-q^{k}}}{1-t^{q^{k}-q}} (t^{q^{k}-q}-1) \right) -
         \prod_{i=1}^{k-1}\dfrac{1-t^{q^{k+1}-q^{i+1}}}{1-t^{q^{k}-q^{i}}} 
         \qbin{n-2}{k-1}{q,t^{q^2}}  \cdot h \\
&=  t^{q^{n}-q^{k}} \qbin{n-2}{k-2}{q,t^{q^2}} -
       \prod_{i=1}^{k-1}\dfrac{1-t^{q^{k+1}-q^{i+1}}}{1-t^{q^{k}-q^{i}}} 
         \qbin{n-2}{k-1}{q,t^{q^2}} \cdot h \\
\end{aligned}
$$
In this last expression, the first summand 
$t^{q^{n}-q^{k}}\qbin{n-2}{k-2}{q,t^{q^2}}$ lies in
$$
(-1)^{(k-2)(n-k)}\NN[t,t^{-1}]=(-1)^{k(n-k)}\NN[t,t^{-1}].
$$
by induction on $n$.

The second summand has three factors, of
which 
\begin{enumerate}
\item[$\bullet$]
the $(q,t)$-binomial lies in $(-1)^{(k-1)(n-k-1)}\NN[t,t^{-1}]$
by induction on $n$, 
\item[$\bullet$]
the product over $i=1,2,\ldots,k-1$ has
each factor in $-\NN[t,t^{-1}]$ as observed in 
\eqref{negative-geometric-series} within
the proof of Case 1, and 
\item[$\bullet$]
the factor of $-h$ also lies in $-\NN[t,t^{-1}]$.
\end{enumerate}
Thus, using the fact that $n-k$ is odd, the second summand lies in 
$$
(-1)^{(k-1)(n-k-1)} \cdot (-1)^{k-1} \cdot (-1)^1\NN[t,t^{-1}]=(-1)^{k(n-k)}\NN[t,t^{-1}].
$$
\end{proof}

\section{Nondegenerate unitary subspaces 
and Ennola duality}
\label{Ennola-section}

There is a well-developed analogy between 
binomial coefficients counting subsets and $q$-binomials counting
subspaces, as well as between the {\it symmetric group} $\symm_n$ and 
{\it general linear group} $GL_n:=GL_n(\FF_q)$.
The authors thank John Shareshian for pointing out an
extension of this analogy to nondegenerate subspaces of $(\FF_{q^2})^n$ as 
a {\it unitary space}\footnote{Meaning that
one equips $(\FF_{q^2})^n$ with a nondegenerate
Hermitian form $(\cdot,\cdot)$, sesquilinear with respect to the conjugation
action $\bar{\alpha}:=\alpha^q$ in $\mathrm{Gal}(\FF_{q^2}/\FF_q)$,
such as $(x,y):=\sum_{i=1}^n x \overline{y}$. See, e.g., Grove \cite[Chapter 10]{Grove} and
Section~\ref{CSP-section} below.},
and the {\it unitary group} $U_n:=U_n(\FF_{q^2})$:
\begin{equation}
\label{Shareshian-observation}
\begin{array}{rcl}
\left( \begin{matrix} n \\ k \end{matrix} \right)
&=& \#\{ \text{subsets of cardinality }k\text{ in }\{1,2,\ldots,n\} \}  \\
&=& [ \symm_n : \symm_k \times \symm_{n-k} ]  \\
& &  \\
\qbin{n}{k}{q}
&=& \#\{ \FF_{q}\text{-subspaces of dimension }k\text{ in }(\FF_{q})^n \}  \\
&=&q^{-k(n-k)}[ GL_n: GL_k \times GL_{n-k} ] =[ GL_n : P_{k,n-k} ]\\
& &  \\
\qbin{n}{k}{-q} 
&=&(-q)^{-k(n-k)}
\#\left\{ \begin{matrix} 
           \text{nondegenerate } \\
          \FF_{q^2}\text{-subspaces of dimension }k\text{ in }(\FF_{q^2})^n
          \end{matrix} 
  \right\}\\
&=&(-q)^{-k(n-k)}[U_n : U_k \times U_{n-k} ].
\end{array}
\end{equation}
Here $P_{k,n-k}$ denotes the parabolic subgroup of $GL_n$ 
that stabilizes one particular choice of a $k$-dimensional
$\FF_q$-subspace, with $GL_k \times GL_{n-k}$ its Levi subgroup of
index $[P_{k,n-k}:GL_k \times GL_{n-k}]=q^{k(n-k)}$.
The formulas in  \eqref{Shareshian-observation}
follow easily from the transitivity of the actions 
on subsets, subspaces, nondegenerate subspaces of the groups $\symm_n, GL_n, U_n$, 
along with these well-known cardinalities (see, for example, 
Grove \cite[Chapters 1, 10, 11]{Grove}):
$$
\begin{aligned}
|\symm_n|&=n!\\
|GL_n|&=q^{\binom{n}{2}} (q-1)(q^2-1)(q^3-1)(q^4-1)\cdots (q^n-1)\\
|U_n|&=q^{\binom{n}{2}} (q+1)(q^2-1)(q^3+1)(q^4-1)\cdots (q^n-(-1)^n).\\
\end{aligned}
$$
Instead, we would like to place \eqref{Shareshian-observation} within
the context of {\it Ennola duality}, relating
unipotent characters of $GL_n$ to those
of the finite {\it unitary group} $U_n:=U_n(\FF_{q^2})$.
We review a portion of this material here--
see Thiem and Vinroot \cite{ThiemVinroot} for a more extensive treatment.

We first review the notion of a unipotent character for $GL_n$ or
$U_n$.  A {\it torus} $T$ in a finite group of Lie type $G$
is an abelian subgroup containing only semisimple elements.  
Given any linear character $\theta:T\rightarrow \CC$
of a {\it maximal} torus $T$, 
there is a virtual character $R_T^G(\theta)$ of $G$ called 
the \emph{Deligne--Lusztig character} of the pair $(T,\theta)$ 
(see \cite[Chapter 11]{DigneMichel}, \cite[Chapter 7]{Carter}).  
A natural $\CC$-subspace of the space of class functions of $G$ is 
$$
\mathcal{U}(G)=\mathbb{\CC}\text{-span}\{R_T^G(1)\mid T\text{ a maximal torus}\},
$$
where $1$ is the trivial character of $T$.
In the case where $G$ is $GL_n$ or $U_n$, for each 
partition $\lambda$ of $n$, there is a unique maximal
torus $T_\lambda$ up to conjugacy, and $\mathcal{U}(G)$  
has a $\CC$-basis given by irreducible
characters which we will denote $\chi_{GL}^\lambda$ and  $\chi_U^\lambda$;
one calls characters in this space {\it unipotent characters}.
We conform here to Macdonald's convention \cite[Chap. IV]{Macdonald}
that $\chi_{GL}^{(1^n)}$ is the trivial representation of $GL_n$;
this differs from some conventions by the 
{\it conjugation} operation
$\lambda \leftrightarrow \lambda'$ on partitions, 
that is, transposing their Ferrers diagrams.

With this convention, if one defines
$$
n(\lambda):=\sum_{i \geq 1}(i-1)\lambda_i,
$$
then the degree of the $GL_n$-character $\chi_{GL}^\lambda$
is a polynomial in $q$, of degree $\binom{n}{2}-n(\lambda)$,
having the following explicit product expression \cite[Chap. IV, (6.7)]{Macdonald}:
\begin{equation}
\label{q-hook-formula}
\chi^\lambda_{GL}(1)=f^\lambda(q):=
q^{n(\lambda')}(q)_n \prod_{(i,j)} \left(1-q^{h_{ij}}\right)^{-1}.
\end{equation}
Here the product runs over $(i,j)$ with $i \geq 1$ and
$j \leq \lambda_i$, that is, over the cells in the Ferrers
diagram for $\lambda$, and $h_{ij}:=\lambda_i-i+\lambda'_j-j+1$ 
is the {\it hooklength} at cell $(i,j)$,
where $\lambda_j'$ is the length of the $j^{th}$ column 
in the diagram.  
The degrees of $\chi_\symm^\lambda$ and $\chi_U^\lambda$
are then determined by the same polynomial $f^\lambda(q)$:  
one has $\chi_\symm^\lambda(1)=f^\lambda(1)$, and 
\begin{equation}
\label{Ennola-duality}
\chi_U^\lambda(1) = (-1)^{\binom{n}{2}-n(\lambda)} f^\lambda(-q).
\end{equation}

The relation \eqref{Ennola-duality} is what we are calling here
{\it Ennola duality}.  We wish to extend it to 
explain \eqref{Shareshian-observation},
utilizing the {\it characteristic maps} isomorphisms
for the three families of groups $G_n=\symm_n, GL_n, U_n$;
for $\symm_n$ see Macdonald \cite[Chap I \S 7]{Macdonald},
for $GL_n$ see Macdonald \cite[Chap IV \S 4]{Macdonald},
and for $U_n$ see and Thiem and Vinroot \cite[\S 4]{ThiemVinroot}.
In each case, these are $\CC$-linear isomorphisms 
$
\UUU(G_n) \overset{\ch_G}{\longrightarrow} \Lambda_n
$ 
where $\Lambda_n$ denotes the space of  
{\it symmetric functions} with $\CC$ coefficients 
which are homogeneous of degree $n$.
The characteristic maps are defined by
$$
\begin{array}{rcl}
 \chi_\symm^\lambda &  \overset{\ch_\symm}{\longmapsto} & s_\lambda \\
\chi_{GL}^\lambda &  \overset{\ch_{GL}}{\longmapsto} & s_\lambda \\
\chi_U^\lambda &  \overset{\ch_{U}}{\longmapsto} & (-1)^{\lfloor \frac{n}{2}\rfloor+n(\lambda)} 
                               s_\lambda. \\
\end{array}
$$
where $s_\lambda$ is the {\it Schur function} indexed by the
partition $\lambda$ of $n$; see \cite[Chap. I \S 3]{Macdonald}. 
It is also worth mentioning that if $p_\lambda$ 
is the {\it power sum} symmetric function corresponding to the partition $\lambda$ of $n$, 
then
$$
R_{T_\lambda}^{G}(1) 
\overset{\ch_{G}}{\longmapsto}
(-1)^{|\lambda|-\ell(\lambda)} p_\lambda \quad \text{ for }G=GL_n\text{ or }U_n. 
$$
From the characteristic map one deduces the 
following extension of \eqref{Ennola-duality}.
\begin{proposition}
\label{extended-Ennola-duality}
Given three class functions $\chi_\symm, \chi_{GL}, \chi_U$
in $\UUU(G_n)$ for the three families $G_n$ above,
whenever they have the {\bf same} symmetric function image 
$$
 \ch_\symm \chi_\symm =  \ch_{GL} \chi_{GL} = \ch_U \chi_U.
$$
then the polynomial $f(q)$ giving 
the degree $\chi_{GL}(1)$ satisfies $\chi_\symm(1) = f(1)$ and
$$
\chi_U(1) = \pm f(-q).
$$
\end{proposition}
\begin{proof}
Expand the symmetric function as $\sum_\lambda c_\lambda s_\lambda$
for some integers $c_\lambda$, and apply the inverse of the
characteristic map isomorphism to get these virtual character
expansions and degrees:
$$
\begin{array}{rcl}
\chi_{GL}&=& \sum_\lambda c_\lambda \chi_{GL}^\lambda \\
f(q):=\chi_{GL}(1)&=& \sum_\lambda c_\lambda f^\lambda(q)\\
 & & \\
\chi_\symm&=& \sum_\lambda c_\lambda \chi_\symm^\lambda \\
\chi_\symm(1)&=& \sum_\lambda c_\lambda f^\lambda(1) = f(1)\\
 & & \\
\chi_U&=& \sum_\lambda c_\lambda (-1)^{\lfloor \frac{n}{2}\rfloor+n(\lambda)}\chi_U^\lambda \\
\chi_U(1)&=&\sum_\lambda c_\lambda 
   (-1)^{\lfloor \frac{n}{2}\rfloor+n(\lambda)} \cdot
    (-1)^{\binom{n}{2}-n(\lambda)} f^\lambda(-q)\\
&=& (-1)^{\lfloor \frac{n}{2}\rfloor+\binom{n}{2}} f(-q).
\end{array}
$$
\end{proof}
To explain \eqref{Shareshian-observation},
we need one further fundamental fact 
(see \cite[Chap. I (7.3), Chap. IV (4.1)]{Macdonald}, \cite[Cor. 4.1]{ThiemVinroot})
about the characteristic maps $\ch_G$ for all $n$:
taken together, they give a {\it ring} (and even Hopf algebra) isomorphism 
$$
\UUU(G):=\bigoplus_{n \geq 0} \UUU(G_n) 
\quad \overset{\ch_G}{\longrightarrow} \quad
\bigoplus_{n \geq 0}\Lambda_n=:\Lambda.
$$
Here the ring of symmetric functions $\Lambda$
is given its usual product,
and the $\CC$-vector space $\UUU(G)$ is endowed with product
structure 
$$
\begin{array}{rcl}
\UUU(G_a) \otimes \UUU(G_b) & \longrightarrow & \UUU_{a+b} \\
\chi_a \otimes \chi_b &\longmapsto 
    & R^{G_{a+b}}_{G_a \times G_b} \left( \chi_a \otimes \chi_b \right)
\end{array}
$$
where $R^{G_{a+b}}_{G_a \times G_b}(-)$ is interpreted in the three cases as 
\begin{enumerate}
\item[$\bullet$]
{\it induction} of characters from $\symm_a \times \symm_b$
to $\symm_{a+b}$,
\item[$\bullet$]
{\it Harish-Chandra induction} of characters 
from $GL_a \times GL_b$ to $GL_{a+b}$, that is,
{\it inflation} from $GL_a \times GL_b$ to $P_{a,b}$ by
composing with the quotient map $P_{a,b} \rightarrow GL_a \times GL_b$,
followed by usual induction from $P_{a,b}$ to $GL_{a+b}$, and
\item[$\bullet$]
{\it Deligne-Lusztig induction} of characters 
from $U_a \times U_b$ to $U_{a+b}$; see, for example, \cite[Chapter 11]{DigneMichel} 
and \cite[Chapter 7]{Carter} for a precise definition.
\end{enumerate}
\noindent
A key point is that these three induction operations multiply
the character degree $(\chi_a \otimes \chi_b)(1)=\chi_a(1) \cdot \chi_b(1)$
by a predictable factor, equal to the right side of \eqref{Shareshian-observation}:
\begin{enumerate}  
\item[$\bullet$]
Induction from $\symm_a \times \symm_b$ to $\symm_{a+b}$ multiplies
degrees by $[\symm_n:\symm_a \times \symm_b]$.
\item[$\bullet$]
Inflation from $GL_a \times GL_b$ to $P_{a,b}$ does
not change degree, while induction from $P_{a,b}$ to $GL_{a+b}$
multiplies degrees by $[GL_{a+b}:P_{a,b}]$.
\item[$\bullet$]
Deligne-Lusztig induction from $U_a \times U_b$ to $U_{a+b}$ 
is known \cite[Proposition 12.17]{DigneMichel}
to multiply degrees by $\pm q^{-ab}[U_{a+b}:U_a \times U_b]$.
\end{enumerate}

We now apply Proposition~\ref{extended-Ennola-duality}
to the trivial degree one character 
$
\triv_{G_k} \otimes \triv_{G_{n-k}}
$
so that
$
R^{G_{n}}_{G_k \times G_{n-k}} 
\left( \triv_{G_k} \otimes \triv_{G_{n-k}} \right)
$
has degree given by the right side of \eqref{Shareshian-observation}.
On the other hand, its image under $\ch_G$ is given by 
$$
\ch \triv_{G_k} \cdot \ch \triv_{G_{n-k}}
=\left( \pm s_{1^k} \right) \left( \pm s_{1^{n-k}} \right)
=\pm s_{1^k} s_{1^{n-k}}
$$
for any of the three families, so that
\eqref{Shareshian-observation} becomes 
a special case of Proposition~\ref{extended-Ennola-duality}.

\section{A cyclic sieving phenomenon for nondegenerate subspaces}
\label{CSP-section}

One original motivation for the
$(q,t)$-binomial in \cite{RSW} was its role in an instance 
of the {\it cyclic sieving phenomenon}, which we recall here.
The {\it finite Grassmannian} of all $k$-dimensional $\FF_q$-subspaces
inside $(\FF_q)^n$ carries an interesting action of a cyclic group
$\ZZ/(q^n-1)\ZZ \cong \FF_{q^n}^\times$:  one embeds
$\FF_{q^n}^\times \hookrightarrow GL_n(\FF_q)$ through any
choice of an $\FF_q$-linear isomorphism $\FF_{q^n} \cong \FF_q^n$.
One can then prove \cite[Theorem 9.4]{RSW} that for an element 
$c$ in $\FF_{q^n}^\times$ of multiplicative
order $d$, the number of $k$-dimensional subspaces preserved by $c$ equals
the $(q,t)$-binomial with $t$ evaluated at any primitive
$d^{th}$ root-of-unity.

In light of this result, and the interpretation for
the negative $q$-binomial in terms of {\it nondegenerate}
unitary subspaces given in \eqref{Shareshian-observation},
one might ask for a similar cyclic sieving phenomenon
involving the $(q,t)$-binomial at negative $q$.
Our goal in this section is such a result when $n$ is odd, 
Theorem~\ref{CSP-theorem} below.  
It involves the action of a certain subgroup $C$ of the cyclic
group $\FF_{q^{2n}}^\times$, acting unitarily on $V=\FF_{q^{2n}}$,
and permuting nondegenerate subspaces.

We begin by introducing a compatible family of 
sesquilinear forms on $V=\FF_{q^{2n}}$ for $n$ odd,
that depend upon upon the choice of scalars $\FF_{q^{2m}}$ 
over which one views $V$ as an $\FF_{q^{2m}}$-vector space.

\begin{definition}
Let $q$ be a prime power, and $n$ an odd integer
with
$n \geq 1$.  For each positive divisor $m$ of $n$,
consider $V=\FF_{q^{2n}}$ as an $\FF_{q^{2m}}$-vector space,
and recall that one has the surjective {\it trace map}
$$
\begin{array}{rcl}
\Tr_{\FF_{q^{2n}}/\FF_{q^{2m}}}: \FF_{q^{2n}} 
&\longrightarrow &
\FF_{q^{2m}} \\
\alpha &\longmapsto & 
\alpha + \alpha^{q^{2m}} +\alpha^{q^{4m}} + \cdots + \alpha^{q^{2(n-m)}} 
\end{array}
$$
Use this to define maps
$$
\begin{array}{rcl}
V \times V & \longrightarrow & \FF_{q^{2m}} \\
(\alpha, \beta)& \longmapsto & 
  (\alpha,\beta)_{\FF_{q^{2m}}}:= 
  \Tr_{\FF_{q^{2n}}/\FF_{q^{2m}}} \left(\alpha \cdot \beta^{q^n} \right) \\
\end{array}
$$
Recall also that $\FF_{q^{2m}}$ is a degree two Galois extension of $\FF_{q^m}$,
and hence the nontrivial element $\bar{\alpha}:=\alpha^{q^m}$
of the Galois group defines conjugation $\alpha \mapsto \overline{\alpha}$ on $\FF_{q^{2m}}$.
\end{definition}

\begin{proposition}
Fix the prime power $q$ and the odd integer $n \geq 1$.

Then for each divisor $m$ of $n$, 
the map $(\cdot, \cdot)_{\FF_{q^m}}$ is 
\begin{enumerate}
\item[$\bullet$] $\FF_{q^{2m}}$-Hermitian with
respect to the conjugation on $\FF_{q^{2m}}$, and 
\item[$\bullet$] nondegenerate as an $\FF_{q^{2m}}$-valued form,
\end{enumerate}
thus endowing $V$ with the structure of an 
$\frac{n}{m}$-dimensional unitary space over $\FF_{q^{2m}}$.
\end{proposition}
\begin{proof}
It is straightforward to check that 
$(\cdot, \cdot):=(\cdot, \cdot)_{\FF_{q^{2m}}}$ 
is an additive function of both arguments,
and an $\FF_{q^{2m}}$-linear function of its first argument, that is,
$$
\begin{array}{rlll}
(\alpha+\alpha', \beta)
&=& (\alpha,\beta) +(\alpha', \beta)& \\
(\alpha, \beta+\beta')
&=& (\alpha,\beta) +(\alpha, \beta')& \\
(c \alpha, \beta)
&=&c (\alpha,\beta) &\text{ for }c \in \FF_{q^{2m}}.\\
\end{array}
$$
In checking the other properties that define a Hermitian
form, it is useful to note that
elements of $c$ in $\FF_{q^{2m}}$ satisfy 
$(c^{q^m})^{q^m}=c$, and since $\frac{n}{m}$  is odd, one also has 
$$
c^{q^n} = \underbrace{\left( \cdots 
              \left( \left( c^{q^m}  \right)^{q^m} \right) 
          \cdots \right)^{q^m}}_{\frac{n}{m}\text{ times}} = c^{q^m}.
$$
To check $\FF_{q^{2m}}$-sesquilinearity in the second argument, 
given $c \in \FF_{q^{2m}}$, one has
$$
\begin{array}{rl}
(\alpha, c\beta)
&= \Tr_{\FF_{q^{2n}}/\FF_{q^{2m}}} \left(\alpha \cdot (c\beta)^{q^n} \right) \\
&= \Tr_{\FF_{q^{2n}}/\FF_{q^{2m}}} \left( c^{q^n} \cdot \alpha \cdot \beta^{q^n} \right) \\
&= c^{q^n} \Tr_{\FF_{q^{2n}}/\FF_{q^{2m}}} \left( \alpha \cdot \beta^{q^n} \right) \\
&= c^{q^m} (\alpha, \beta) \\
&= \overline{c} \cdot (\alpha, \beta) 
\end{array}
$$
where the middle equality
used the $\FF_{q^{2m}}$-linearity of $\Tr_{\FF_{q^{2n}}/\FF_{q^{2m}}}$.
One also has
$$
\begin{array}{rll}
\overline{(\beta, \alpha)}
&= \left(
     \Tr_{\FF_{q^{2n}}/\FF_{q^{2m}}} \left(\beta \cdot \alpha^{q^n} \right)
   \right)^{q^m} \\
&=\left(
     \Tr_{\FF_{q^{2n}}/\FF_{q^{2m}}} \left(\beta \cdot \alpha^{q^n} \right)
   \right)^{q^n} \\
&=\Tr_{\FF_{q^{2n}}/\FF_{q^{2m}}} \left(\beta^{q^n} \cdot \left(\alpha^{q^n}\right)^{q^n} \right)\\
&=\Tr_{\FF_{q^{2n}}/\FF_{q^{2m}}} \left(\alpha \cdot \beta^{q^n} \right)\\
&=(\alpha,\beta).
\end{array}
$$
Lastly, one needs to know that $(\cdot,\cdot)$ is nondegenerate
as a pairing on $V$, or equivalently, that for any nonzero
$\alpha$ in $V=\FF_{q^{2n}}$, the $\FF_{q^{2m}}$-linear functional
$$
\begin{array}{rcl}
V &\longrightarrow &\FF_{q^{2m}} \\
\beta &\longmapsto & 
   (\alpha,\beta) =  
    \Tr_{\FF_{q^{2n}}/\FF_{q^{2m}}} \left(\alpha \cdot \beta^{q^n} \right)
\end{array}
$$
is surjective, or equivalently, not identically zero.  
This follows because a separable field extension $K/k$, 
such as $\FF_{q^{2n}}/\FF_{q^{2m}}$, always has nondegenerate
$k$-bilinear pairing 
$
\langle \alpha,\beta \rangle_k := 
\Tr_{K/k}(\alpha \cdot \beta),
$
and $\beta \mapsto \beta^{q^n}$ is
an automorphism of $K=\FF^{q^{2n}}$.
\end{proof}

We note here the following compatibility between the
the various forms $(\cdot, \cdot)_{\FF_{q^{2m}}}$, which will be
used in the proof of Theorem~\ref{CSP-theorem}.

\begin{proposition}
\label{nondegeneracy-consistency}
Fix the prime power $q$ and odd positive integers $\ell,m,n$,
with $\ell$ dividing $m$ and $m$ dividing $n$.

Then an $\FF_{q^{2m}}$-subspace $W$ of $\FF_{q^{2n}}$, 
when regarded as an $\FF_{q^{2\ell}}$-subspace,
is nondegenerate with respect to the form $(\cdot, \cdot)_{\FF_{q^{2m}}}$
if and only if it is nondegenerate with respect to the form 
$(\cdot, \cdot)_{\FF_{q^{2\ell}}}$.
\end{proposition}
\begin{proof}
As in the previous proof, $W$ is 
$(\cdot, \cdot)_{\FF_{q^{2m}}}$-nondegenerate 
if and only if for every nonzero $\alpha$ in $W$, the
$\FF_{q^{2m}}$-linear functional 
$f_{m,\alpha}:W \rightarrow \FF_{q^{2m}}$
given by 
$$
f_{m,\alpha}(\beta)=(\alpha, \beta)_{\FF_{q^{2m}}}
$$
is surjective, or equivalently, not identically zero.
As the corresponding functional 
$f_{\ell,\alpha}:W \rightarrow \FF_{q^{2\ell}}$
can be expressed as the composite
$
f_{\ell,\alpha} = \Tr_{\FF_{q^{2m}}/\FF_{q^{2\ell}}} \circ f_{m,\alpha},
$
where the trace map 
$\Tr_{\FF_{q^{2m}}/\FF_{q^{2\ell}}}: \FF_{q^{2m}} \longrightarrow \FF_{q^{2\ell}}$
is well-known to be surjective,
$f_{m,\alpha}$ is nonzero if and only if $f_{\ell,\alpha}$ is nonzero.
\end{proof}

We next describe the subgroup of the multiplicative
group $\FF_{q^{2n}}^\times$ which will act unitarily with respect
to our chosen Hermitian forms.  Let $\gamma$ be a 
generator for $\FF_{q^{2n}}^\times\cong \ZZ/(q^{2n}-1)\ZZ$ 
as a cyclic group.

\begin{proposition}
Fix $q,n$ as before, and any divisor $m$ of $n$.

Then the power $\gamma^{q^n-1}$ generates a cyclic subgroup 
$C \cong  \ZZ/(q^{n}+1)\ZZ$ of $\FF_{q^{2n}}^\times$
which acts on $V=\FF_{q^{2n}}$ unitarily with respect to
the $\FF_{q^{2m}}$-Hermitian form $(\cdot,\cdot)_{\FF_{q^{2m}}}$.
\end{proposition}
\begin{proof}
The cardinality of 
$C=\langle \gamma^{q^n-1} \rangle$ should be clear.  For the rest, calculate
$$
\begin{aligned}
(\gamma^{q^n-1} \alpha, \gamma^{q^n-1} \beta)_{\FF_{q^{2m}}}
&= \Tr_{\FF_{q^{2n}}/\FF_{q^{2m}}} \left(\gamma^{q^n-1} \alpha \cdot (\gamma^{q^n-1} \beta)^{q^n} \right) \\
&= \Tr_{\FF_{q^{2n}}/\FF_{q^{2m}}} \left((\gamma^{q^n-1})^{q^n+1} \cdot \alpha \cdot \beta^{q^n} \right) \\
&= \Tr_{\FF_{q^{2n}}/\FF_{q^{2m}}} \left(\gamma^{q^{2n}-1} \cdot \alpha \cdot \beta^{q^n} \right) \\
&=(\alpha,\beta)_{\FF_{q^{2m}}}.
\end{aligned}
$$
\end{proof}

To state our cyclic sieving phenomenon, fix $n$ odd as above,  
a prime power $q$, and a $k$ in the range $0 \leq k \leq n$.
Consider the set 
$$
X:=
\left\{ 
\begin{matrix}
\text{ all }(\cdot,\cdot)_{\FF_{q^2}}\text{-nondegenerate }k\text{-dimensional }\ \\
\FF_{q^2}\text{-subspaces of }V=\FF_{q^{2n}} 
\end{matrix}
\right\}.
$$
Since the group 
$$
C:=\langle \gamma^{q^n-1} \rangle \cong \ZZ/(q^n+1)\ZZ
$$ 
acts unitarily with respect to this form, $C$ permutes the set $X$.

To define a polynomial $X(t)$ in $\NN[t]$ we first define a 
$t$-version of $q^{k(n-k)}\qbin{n}{k}{q}$ by
$$
Y(q,t):=\left( \prod_{i=0}^{k-1}\frac{1-t^{q^n-q^{n-k+i}}}{1-t^{q^k-q^{i}}} \right)
\qbin{n}{k}{q,t}.
$$
We define $X(t)$ as a polynomial version of $Y(-q,t)$, namely
$$
X(t):=t^E \cdot 
\prod_{i=0}^{k-1}\frac{1-t^{q^n+(-q)^{n-k+i}}}{1-t^{q^k-(-1)^k(-q)^{i}}}
\cdot 
\prod_{i=0}^{k-1}\frac{1-t^{q^n+(-q)^{i}}}{1-t^{q^k-(-1)^k(-q)^{i}}},
$$
where one defines
$$
E:=\begin{cases}
0 &{\text{ if $k$ is odd,}}\\
2\sum_{i=0}^{k-1} (q^k-(-q)^i)& {\text{ if $k$ is even.}}
\end{cases}
$$
One may show that Theorem \ref{lastthm} implies $X(t)$ lies in $\NN[t]$. 
It is also easily checked that, since $n$ is odd, one has
\begin{equation}
\label{negative-(q,t)-bin-evaluates-right-at-1}
X(1)= (-q)^{k(n-k)}\qbin{n}{k}{-q}.
\end{equation}
Furthermore, whenever $q$ is odd, all powers of $t$ in
$X(t)$ occur with even exponents, and thus one has
\begin{equation}
\label{negative-(q,t)-bin-evaluates-right-at-(-1)}
X(-1)=X(1) \text{ for odd }q.
\end{equation}

\begin{theorem}
\label{CSP-theorem}
This triple $(X,X(t),C)$ exhibits a {\sf cyclic sieving phenomenon}
as in \cite{RSW}:  for any $c$ in $C$,  
the number of elements $x$ in $X$ having $c(x)=x$ is given by
evaluating $X(t)$ with $t$ any complex root-of-unity whose
multiplicative order is the same as $c$.
\end{theorem}
\begin{proof}
Given $c$ in $C\subset \FF_{q^{2n}}^\times$, we wish to count how many
$x$ in $X$ have $c(x)=x$.  Let $\FF_{q^2}(c)$ denote the subfield of
$\FF_{q^{2n}}$ generated by $\FF_{q^2}$ and $c$, so there exists
a unique divisor $m$ of $n$ for which 
\begin{equation}
\label{m-defining-equation}
\FF_{q^2}(c)=\FF_{q^{2m}}.
\end{equation}
A $k$-dimensional $\FF_{q^2}$-subspace 
$W \subset V=\FF_{q^{2n}}$ is fixed by $c$ if and only if $cW \subset W$, that is,
if and only if $W$ is actually a subspace over $\FF_{q^2}(c)$.
By \eqref{m-defining-equation}, this is equivalent to $W$ being
a $k'$-dimensional $\FF_{q^{2m}}$-subspace, where $k':=\frac{k}{m}$.

According to Proposition~\ref{nondegeneracy-consistency},
this $\FF_{q^2}$-subspace $W$ is in addition nondegenerate for 
$(\cdot,\cdot)_{\FF_{q^2}}$ if and only if it is nondegenerate for
$(\cdot,\cdot)_{\FF_{q^{2m}}}$.
Therefore the number of $x$ in $X$ with $c(x)=x$ will be the number of
 $(\cdot,\cdot)_{\FF_{q^{2m}}}$-nondegenerate $k'$-dimensional
$\FF_{q^{2m}}$-subspaces of $V=\FF_{q^{2n}} \cong (\FF_{q^{2m}})^{n'}$
where $n':=\frac{n}{m}$.  By \eqref{Shareshian-observation}, this number is 
\begin{equation}
\label{predicted-CSP-value}
(-Q)^{k'(n'-k')}
\qbin{n'}{k'}{-Q}, \quad \text{ where }Q:=q^m.
\end{equation}

On the other hand, assuming that $c$ has multiplicative order $A$, one
can evaluate $X(t)$ at $t=\omega$ a primitive $A^{th}$ root-of-unity. This 
makes heavy use of Proposition~\ref{numbth} below, relating $A$
to the number $m$ defined by  \eqref{m-defining-equation} above, and
allowing one to analyze the locations of zeroes in the
numerator and denominators appearing in the explicit formula for 
$X(t)$.  In particular, it will be shown that
$X(\omega)$ vanishes unless $m$ divides both $n$ and $k$,
and then a limiting procedure will yield the predicted 
value \eqref{predicted-CSP-value} for $X(\omega)$ in this case.  
We proceed in cases based on the value of $A$.

\vskip.1in
\noindent
{\sf Case 1.  $A=1$}.

This case follows from equality \eqref{negative-(q,t)-bin-evaluates-right-at-1}
combined with \eqref{Shareshian-observation}.

\vskip.1in
\noindent
{\sf Case 2.  $A=2$}.

If $A=2$ then $c=-1 \neq +1$ in $\FF_{q^{2n}}$, forcing $q$ to be odd.
In this case, $cW=W$ for all subspaces $W$, and 
\eqref{negative-(q,t)-bin-evaluates-right-at-1} shows that
$X(-1)=X(1)$, so the result follows as in the $A=1$ case.

\vskip.1in
\noindent
{\sf Case 3.  $A \geq 3$}.

Given $A$, let $m$ be as in \eqref{m-defining-equation}.
We first show $X(\omega)=0$ if $m$ does not divide $k$.
Note that the first product in the definition of $X(t)$, namely
\begin{equation}
\label{first-product-in-X(t)}
\prod_{i=0}^{k-1}\frac{1-t^{q^n+(-q)^{n-k+i}}}{1-t^{q^k-(-1)^k(-q)^{i}}}
\end{equation}
has each of its factors a polynomial in $t$, 
since $n$ being odd implies
$$
q^n+(-q)^{n-k+i}= q^{n-k} ( q^k-(-1)^k(-q)^{i} ).
$$
Thus \eqref{first-product-in-X(t)} never has poles.
As for the second product in the definition of $X(t)$, namely
\begin{equation}
\label{second-product-in-X(t)}
\prod_{i=0}^{k-1}\frac{1-t^{q^n+(-q)^{i}}}{1-t^{q^k-(-1)^k(-q)^{i}}},
\end{equation}
it has powers of $t$ with exponents
$$
\begin{array}{ll}
q^{n}+1, \,\, q^n-q, \,\, q^n-q^2, \ldots, \,\, q^n+(-1)^{k-1} q^{k-1} 
&\text{ in the numerator},\\
q^k+1, \,\, q^k-q, \,\, q^k+q^2, \,\, \ldots, \,\, q^k+q^{k-1} &\text{ in the denominator for odd }k, \\
q^k-1, \,\, q^k+q, \,\, q^k-q^2, \,\, \ldots, \,\, q^k+q^{k-1} &\text{ in the denominator for even }k.
\end{array}
$$
Since $n$ is an odd multiple of $m$, 
Proposition \ref{numbth}(iii,iv) implies that the choice $t=\omega$ yields 
$\lceil k/m \rceil$ numerator zeroes, from exponents
$$
q^n+1, q^n-q^m, q^n+q^{2m}, q^n-q^{3m}, \ldots
$$
and $\lfloor k/m \rfloor$ denominator zeroes, from exponents
$$
q^k+q^{k-m}, q^k-q^{k-2m}, q^k+q^{k-3m}, \ldots
$$
regardless of the parity of $k$.
Thus when $m$ does not divide $k$, the numerator has more zeros than the 
denominator, and $X(\omega)=0$.

When $m$ does divide $k$, we wish to evaluate $X(t)$
at $t=\omega$ a primitive $A^{th}$ root-of-unity, 
using this general L'H\^opital's rule calculation: if $r \equiv \pm s \bmod{A}$ then
\begin{equation}
\label{L'Hopital-calculation}
\lim_{t \rightarrow \omega} \frac{1-t^r}{1-t^s} =
\begin{cases}
r/s & \text{ if }r \equiv s \equiv 0  \bmod{A}\\
1           & \text{ if }r \equiv s \not\equiv 0  \bmod{A}\\
-\omega^{-s} & \text{ if }r \equiv -s \not\equiv 0  \bmod{A}.\\
\end{cases}
\end{equation}
Pairing zeroes at $t=\omega$ in numerator, denominator
of \eqref{second-product-in-X(t)} and 
using \eqref{L'Hopital-calculation} yields
$$
\prod_{i=0}^{k/m-1} \frac{q^n-(-q^m)^i}{q^k-q^k(-q^{-m})^{i+1}}=
\qbin{n'}{k'}{-Q}.
$$

One can do a similar analysis for the first factor  \eqref{first-product-in-X(t)}
evaluated at $t=\omega$.  This time one finds exponents on $t$ of
$$
\begin{array}{llllll}
q^{n}+q^{n-1}, & q^n-q^{n-2}, & q^n+q^{n-3}, &\cdots,  & q^n+q^{n-k}
&\text{ in the numerator, and}\\
q^k+1,       & q^k-q,      & q^k+q^2,  & \cdots, & q^k+q^{k-1}
&\text{ in the denominator.}
\end{array}
$$
The corresponding zeroes are 
$$
\begin{aligned}
q^n+q^{n-m}, q^n-q^{n-2m}, q^n+q^{n-3m}, \cdots
&\text{ in the numerator, and} \\
q^k+q^{k-m}, q^k-q^{k-2m}, q^k+q^{k-3m}, \cdots
&\text{ in the denominator}
\end{aligned}
$$
whose limit using \eqref{L'Hopital-calculation} yields.
$$
\prod_{i=0}^{k/m-1} \frac{q^n-q^n(-q^{-m})^{i+1}}
{q^k-q^k(-q^{-m})^{i+1}}=q^{(n-k)k/m}=Q^{k'(n'-k')}.
$$

It only remains to analyze the {\it nonzero} factors at $t=\omega$
in the numerator and denominators of the two products
comprising $X(t)$.  We treat this in two cases based on the parity of $k$.

For $k$ odd, we claim that these nonzero numerator and denominator
factors in the second product \eqref{second-product-in-X(t)} pair off 
to give $1$ using \eqref{L'Hopital-calculation}.  To see this claim,
note that since $A$ divides $q^m+1$, one has
$\omega^{q^m+1}=1$, and so one needs only check that the difference of the 
$t$-exponents 
$$
q^n+(-q)^i-(q^k+(-q)^i)=q^k(q^{n-k}-1)
$$ 
is divisible by $q^m+1$. But $m$, which is odd, 
divides $n-k$, which is even, so one also has $2m$ dividing $n-k$. 
Thus $q^{2m}-1$ divides $q^{n-k}-1$, as does $q^m+1$.
We similarly claim that, for $k$ odd, 
the nonzero numerator and denominator factors in the first product
\eqref{first-product-in-X(t)} pair off to give factors
of $1$ using \eqref{L'Hopital-calculation}, because 
the difference of the exponents 
$$
q^n+(-q)^{n-k+i}-(q^k+(-q)^i)=(q^k+(-q)^i)(q^{n-k}-1)
$$ 
is again divisible by $q^m+1$.

For $k$ even, we claim that the nonzero numerator and denominator
factors in the second product \eqref{second-product-in-X(t)} pair off 
in such a way that one can apply the third case of \eqref{L'Hopital-calculation}:
one has as sum of numerator and denominator $t$-exponents
$$
\left( q^n+(-q)^i\right) + \left( q^k-(-q)^i\right)
=q^k(q^{n-k}+1) \equiv 0 \bmod{q^m+1}
$$ 
where the congruence follows as $m$ divides $n-k$ and both are odd. 
Each such factor contributes $-\omega^{-(q^k-(-q)^i)}$ by \eqref{L'Hopital-calculation},
and there are $k-\frac{k}{m}$ such factors, an even number
since $k$ is odd and $m$ is even, giving a total contribution of 
$\omega^{-\sum_{i=0}^{k-1} (q^k-(-q)^i)}$.
Similarly, when $k$ is even,  we claim that these nonzero numerator and denominator
factors in the first product \eqref{first-product-in-X(t)} pair off 
with the sum of the numerator and denominator $t$-exponents 
$$
\left( q^n+(-q)^{n-k+i}\right) + (q^k+(-q)^i)=(q^k-(-q)^i)(q^{n-k}+1) \equiv 0 \bmod{q^m+1}
$$
where the congruence follows for the same reason.  Again there are
$k-\frac{k}{m}$ such factors, giving a total contribution of
$\omega^{-\sum_{i=0}^{k-1} (q^k-(-q)^{i})}$. 

Together these contribute $\omega^{-2\sum_{i=0}^{k-1} (q^k-(-q)^{i})}$,
cancelled by $t^E$ for $k$ even.
\end{proof} 

The following proposition collects technical facts used in
the preceding proof.

\begin{proposition}
\label{numbth} Assume that $c$ in $\FF_{q^{2n}}^\times$
has multiplicative order $A$ at least $3$, and that
$\FF_{q^2}(c)=\FF_{q^{2m}}$, where $m$ divides $n$.
\begin{enumerate}
\item[(i)] The order $A$ must divide $q^m+1$. 
\item[(ii)] The smallest positive integer $d$ such that $A$ 
divides $q^d+1$ is $m.$
\item[(iii)] The order $A$ divides $q^s+q^t$ if and only if 
$s-t$ is an odd multiple of $m$.
\item[(iv)] The order $A$ divides $q^s-q^t$ if and only if 
$s-t$ is an even multiple of $m$.
\end{enumerate}
\end{proposition}
\begin{proof}
{\sf Assertions (i) and (ii).}
These will be deduced from the stronger
\begin{quote}
{\bf Claim}: 
If $d|n$ and $A|q^n+1$, then $A|q^{2d}-1$  if and only if $A|q^d+1$.
\end{quote}
The ``if'' direction in the claim
is clear, so we must only show that $A|q^{2d}-1$ implies $A|q^d+1$. 
If $d=n$ this is the hypothesis on $A$.  So assume that $d<n$, and 
since $n$ is odd, we have $2d<n$.
Clearly $A$ must divide $\gcd(q^n+1,q^{2d}-1)$, which we now 
prove is $q^d+1$. Expressing
$$
q^n+1=q^{n-2d}(q^{2d}-1)+q^{n-2d}+1
$$ 
in order to use Euclidean algorithm, one has
$$
\gcd(q^n+1,q^{2d}-1)=\gcd(q^{n-2d}+1,q^{2d}-1)
$$
Since $d$ divides $n$, and $n$ is odd we have $n\equiv d\bmod{2d}$, and 
therefore 
$$
\gcd(q^n+1,q^{2d}-1)=\gcd(q^{d}+1,q^{2d}-1)=q^d+1.
$$
Given the claim, assertions (i),(ii) follow, since
$\FF_{q^2}(c)=\FF_{q^{2m}}$ means that $m$ is the smallest positive integer $d$ such that 
$c\in \FF_{q^{2d}}$, i.e., such that $A$ divides $q^{2d}-1$. 

\vskip.1in
\noindent
{\sf Assertion (iii).}
In one direction, if $s-t$ is an odd multiple of $m$, 
then $q^m+1$ divides $q^{s-t}+1$ and also
$q^{s}+q^t$, so $A$ also divides  $q^{s}+q^t$.  

For the converse, suppose that $A$ divides  
$q^{s}+q^t$, and assume without loss of generality that
$s \geq t$.  Since $\gcd(A,q)=1$, one has
that $A$ also divides $q^{s-t}+1$ and $q^m+1$,
so $A$ divides $\gcd(q^{s-t}+1,q^m+1)$.
Since $A \geq 3$ we can assume that $s>t$
and write $s-t= m\alpha+\beta$ with $0\le \beta<m$. Expressing
$$
q^{s-t}+1=(q^{m(\alpha-1)+\beta}-q^{m(\alpha-2)+\beta}+\cdots+ 
(-1)^{\alpha-1}q^\beta)(q^m+1)+(-1)^{\alpha}q^\beta+1
$$
and using the Euclidean algorithm, one concludes that
$A$ divides $(-1)^\alpha q^\beta+1$. 

If $\alpha$ is even, then $0<\beta<m$ 
contradicts the minimality of $m$, while $\beta=0$ contradicts $A \geq 3$. 
Thus $\alpha$ is odd, and $A$ divides $q^\beta-1$. 

Now if $\beta=0$, then $s-t$ is an odd multiple of $m$, so we are done. 
Otherwise, if $\beta>0$, then write $m=\gamma\beta+\delta$ 
with $0\le \delta<\beta<m$. 
Expressing  
$$
q^{m}+1=(q^{\beta(\gamma-1)+\delta}+q^{\beta(\gamma-2)+\delta}+\cdots+ 
q^\delta)(q^\beta-1)+q^\delta+1
$$
and using the Euclidean algorithm, one concludes that
$A$ divides $q^\delta+1$.
Again by minimality of $m$ this implies that $\delta=0$, which 
contradicts $A \geq 3$.

\vskip.1in
\noindent
{\sf Assertion (iv).}
In one direction, if $s-t$ is an even multiple of $m$, 
then $q^m+1$ divides $q^{2m}-1$, and hence also
divides $q^{s-t}-1$, and therefore divides $q^{s}-q^t$. 

For the converse, suppose that $A$ divides  
$q^{s}-q^t$, and assume without loss of generality that $s \geq t$.
Since $\gcd(A,q)=1$, one has that $A$ also divides $q^{s-t}-1$ and $q^m+1$,
so $A$ divides $\gcd(q^{s-t}-1,q^m+1).$
Again writing $s-t= m\alpha+\beta$ with $0\le \beta<m$, and
expressing
$$
q^{s-t}-1=(q^{m(\alpha-1)+\beta}-q^{m(\alpha-2)+\beta}+\cdots+ 
(-1)^{\alpha-1}q^\beta)(q^m+1)+(-1)^{\alpha}q^\beta-1
$$
the Euclidean algorithm implies that 
$A$ divides $(-1)^\alpha q^\beta-1$.

If $\alpha$ is odd, this contradicts the minimality $m$ for $\beta>0$, so 
we can assume $\alpha$ is even. The argument proceeds as for Assertion (iii), 
$\beta=0$, and $s-t$ is an even multiple of $m$.
\end{proof}

\section{Remarks and further questions}
\label{remarks-section}

\subsection{Reformulating Theorem~\ref{mainth}
via partitions}
It is well-known (see e.g. \cite[p. 40]{A})
that $\Omega_{n,k}$ bijects with integer partitions $\lambda$ 
whose Ferrers diagram lie inside an $(n-k) \times k$ rectangle.
One version of this bijection sends the word $\omega=(\omega_1,\ldots,\omega_n)$
to the partition $\lambda$
whose Ferrers diagram (drawn in the plane
$\ZZ^2$ in English notation) has its northwest corner at 
$(0,n-k)$, and whose outer boundary is the lattice path 
from $(n-k,k)$ to $(0,0)$ having its $i^{th}$ step go one down (resp.
leftward) if $\omega_i=0$ (resp. $\omega_i=1$).  One has 
$\inv(\omega)=|\lambda|=\sum_i \lambda_i$, the weight of $\lambda$.

We omit the details in verifying the following.
\begin{proposition}
 \label{mainthpart}
Under the above bijection, one has the following
correspondences.
\begin{enumerate}
\item[(i)]
The subset $\Omega_{n,k}' \subset \Omega_{n,k}$ corresponds
to those partitions $\lambda$ inside $(n-k)\times k$ for which 
 \begin{enumerate}
 \item if $k$ is even, each odd part has even multiplicity,
 \item if $k$ is odd, each even part has even multiplicity, and moreover the
 number of parts has the same parity as $n-k$.
 \end{enumerate}
\item[(ii)]
The statistic $p(\omega)$ counting
occurrences of paired $\underline{10}$ in $\omega$ 
corresponds to the statistic $p(\lambda)$ counting
the corner cells in $\lambda$ that are {\bf special}
in the following sense:  they are the last cells in rows
of $\lambda$ corresponding to the last occurrences
of each part with the same parity as $k$.
\item[(iii)]
Theorem~\ref{mainth} becomes
$$
\qbin{n}{k}{q}' = 
  \sum_{ \lambda } q^{|\lambda|-p(\lambda)} (q-1)^{p(\lambda)}.
$$
where the sum runs over $\lambda \subset (n-k) \times k$
satisfying condition (i) above.
\end{enumerate}
\end{proposition}

\noindent
Here are three examples of assertion (iii), with the
first compared to the example appearing just
after Theorem~\ref{mainth}:
$$
\begin{matrix}
(n,k)=(5,2)\\
\begin{tabular}{|l|l|l|}\hline
$\omega \in \Omega'_{5,2}$&$\lambda$ & $\wt(\lambda)$ \\ \hline\hline
$\underline{1}\ \underline{10}\ \underline{0}\ \underline{0}$&
$222$ & $q(q-1)q^2q^2$\\\hline
$\underline{0}\ \underline{1}\ \underline{10}\ \underline{0}$&
$22$ & $q(q-1)q^2$ \\ \hline
$\underline{1}\ \underline{00}\ \underline{10}$
&$211$ & $q(q-1)q^2$ \\ \hline 
$\underline{0}\ \underline{0}\ \underline{1}\ \underline{10}$ &
$2$ & $q(q-1)$\\ \hline
$\underline{0}\ \underline{1}\ \underline{00}\ \underline{1}$ &
$11$ & $q^2$ \\ \hline
$\underline{0}\ \underline{0}\ \underline{0}\ \underline{1}\ \underline{1}$ &
$\emptyset$ & $1$\\ \hline
\end{tabular}
\end{matrix}
\,\,\,
\begin{matrix}
(n,k)=(5,3)\\
\begin{tabular}{|l|l|}\hline
$\lambda$ & $\wt(\lambda)$ \\ \hline\hline
$33$& $q^2(q-1)q^3$ \\\hline
$22$& $q^2q^2$ \\\hline
$31$& $q^2(q-1)(q-1)$ \\\hline
$11$& $(q-1)q^1$ \\\hline
$\emptyset$ & $1$ \\\hline
\end{tabular}
\end{matrix}
\,\,\,
\begin{matrix}
(n,k)=(6,3)\\
\begin{tabular}{|l|l|}\hline
$\lambda$ & $\wt(\lambda)$ \\ \hline\hline
$333$& $q^2(q-1)q^3q^3$ \\\hline
$322$& $q^2(q-1)q^2q^2$ \\\hline
$331$& $q^2(q-1)q^3(q-1)$ \\ \hline
$221$& $q^2q^2(q-1)$\\\hline
$311$& $q^2(q-1)(q-1)q$ \\\hline
$3$&   $q^2(q-1)$ \\\hline
$111$& $(q-1)q^2$ \\ \hline
$1$&   $q-1$\\\hline
\end{tabular}
\end{matrix}
$$

\subsection{Reformulating Theorem~\ref{mainth}
via subspaces}
When $q$ is a prime power, so the size of the finite field $\FF_q$,
one can also reformulate Theorem~\ref{mainth}
as counting certain $k$-dimensional $\FF_q$-subspaces of $\FF_q^n$.  

Recall that a $k$-dimensional subspace $V$ is the column-space of a 
matrix $A$ in $\FF_q^{n \times k}$ in {\it column-echelon form}:
\begin{enumerate}
\item[$\bullet$]
each column ends with a string of $0$'s, preceded by a pivot entry $1$, 
\item[$\bullet$]
with only zeroes in the same row as any pivot, and 
\item[$\bullet$]
where the row indices of the pivots decreasing
from left-to-right.
\end{enumerate}
The map $f$ sending $V$ to the word
$\omega$ in $\Omega_{n,k}$ whose ones are in the same positions
as the row indices of the pivots of $A$ corresponds (see e.g.
\cite{NSW}, \cite[Chapter 1]{Stanley}) to the Schubert
cell decomposition of the Grassmannian $\GG(k,\FF_q^n)$.

If the word $\omega$ corresponds as above to the partition
$\lambda$ inside $(n-k) \times k$, then there
are $q^{\inv(\omega)}=q^{\lambda}$ elements in the pre-image $f^{-1}(\omega)$:
deleting the $k$ pivot rows from $A$ gives an $(n-k) \times k$ matrix
whose nonzero entries lie in the cells of $\lambda$.  
As an example, consider matrices $A$ with this column-echelon form
for $(n,k)=(11,5)$:
$$
\left[
\begin{matrix}
0 & 0 & 0 & 0 & 1 \\
0 & 0 & 0 & 1 & 0 \\
* & * & * & 0 & 0 \\
* & * & * & 0 & 0 \\
* & * & * & 0 & 0 \\
* & * & \underline{*} & 0 & 0 \\
0 & 0 & 1 & 0 & 0 \\
0 & 1 & 0 & 0 & 0 \\
* & 0 & 0 & 0 & 0 \\
\underline{*} & 0 & 0 & 0 & 0 \\
1 & 0 & 0 & 0 & 0 \\
\end{matrix}
\right]
\quad \text{ has }\quad
\lambda =
\begin{matrix}
* & * & * \\
* & * & * \\
* & * & * \\
* & * & \underline{*}\\
* &   &  \\
\underline{*} &  &  
\end{matrix}
\quad \text{ and } \quad
\omega=
\underline{10}\
\underline{0}\
\underline{1}\
\underline{10}\
\underline{0}\
\underline{0}\
\underline{0}\
\underline{1}\
\underline{1}.
$$
The {\it special} entries in 
the matrix $A$, corresponding to the special corner
cells of $\lambda$ from Proposition~\ref{mainthpart}
and corresponding to the paired $\underline{10}$'s in $\omega$, are
shown underlined.

One may then analogously reintepret 
Theorem \ref{mainth} (or Theorem \ref{mainthpart}) 
as saying that the primed $q$-binomial counts
those $k$-dimensional subspaces $V$ whose
column-echelon form has all {\it special entries nonzero}.

Using echelon forms, it was shown in \cite[\S 5.3]{RS} how to associate
to each $k$-dimensional subspace $V$ a power of $t$ so that their
generating function in $t$ interprets the $(q,t)$-binomial coefficient.
Similarly, one can use Theorem~\ref{lastthm} and its proof
to associate a power of $t$ to each such subspace $V$ having 
special entries nonzero,
so as to give a generating function interpretation to the $(q,t)$-binomial
coefficient when $q$ is a negative integer.  
We omit this formulation here.

\subsection{Geometry}
Given a field $\FF$, let $X_\FF$ denote the Grassmannian of
$k$-dimensional subspaces in $\FF^n$.  Aside from its interpretation
when $q$ is a prime power as counting the points of the finite
Grassmannian $X_{\FF_q}$, the $q$-binomial coefficient has 
two well-known interpretations as the Poincar\'e polynomials
$$
\begin{aligned}
&\sum_{i} \rank_\ZZ H^{2i}(X_\CC;\ZZ) q^i \\
&\sum_{i} \dim_{\FF_2} H^{i}(X_\RR;\FF_2) q^i \\
\end{aligned}
$$
See, e.g., \cite{Takeuchi} for the second interpretation.
These lead to the following interpretations
for at least the $q=1$ specialization of the
primed $q$-binomial 
$$
\qbin{n}{k}{q=1}'=
\#\{\omega\in\Omega': p(\omega)=0\}
$$

\begin{enumerate}
\item[$\bullet$]
as the {\it signature} or {\it index} of $X_\CC$ (see \cite{R-Eng}), and
\item[$\bullet$]
as the {\it Euler characteristic} of $X_\RR$, up to a $\pm$ sign.
\end{enumerate}

\begin{question}
\label{q}
Can one generalize either of
the above geometric interpretations for its $q=1$ specialization
to a geometric interpretation for the full primed $q$-binomial?
\end{question}

\subsection{Lack of monotonicity for $(q,t)$}
As mentioned earlier, Theorem~\ref{mainth} makes
inequality \eqref{ineq} transparent.  Thus one
might hope for an analogous inequality involving the
$(q,t)$-binomial and its $-q$ relative, perhaps via
Theorem~\ref{lastthm}.  

Unfortunately, a naive guess along these lines fails 
already for $(n,k,q)=(4,2,4)$:  even though
$\qbin{n}{k}{q,t}$ and $(-1)^{k(n-k)} \qbin{n}{k}{-q,t}$ both
lie in $\NN[t]$ and have the same degree $k(q^n-q^k)=480$, 
their difference 
contains both positive and negative coefficients.

\subsection{A conjecture on Schur functions}

We conjecture a generalization of Theorem~\ref{lastthm} that applies
to a $(q,t)$-analogue $S_\lambda(1,t,\ldots,t^n)$
of principally specialized Schur functions, discussed
in \cite[Definition 5.1]{RS}.  Define for
integer partitions $\lambda$, the statistic
$b(\lambda):= \sum_{i} (i-1)\lambda_i.$

\begin{conjecture}
\label{conj}
If $q\le -2$ is a negative integer, then
$$
(-1)^{n|\lambda|-b(\lambda)}S_{\lambda}(1,t,\ldots,t^n)
$$
is a Laurent polynomial in $t$,
all of whose coefficients are non-negative integers.
\end{conjecture}

\subsection{Generating function for $a,p$ on $\Omega'_{n,k}$}
One may find an explicit product representation for the rational 
generating function
$$
G_k(x,q,z)=\sum_{n\ge k}x^n \sum_{\omega\in \Omega'_{n,k}} 
q^{a(\omega)} z^{p(\omega)}.
$$
We do not give the result here, but note one of its specializations
$$
G_k(x,1,1)=\sum_{n\ge k}x^n  |\Omega'_{n,k}|=
\frac{x^k}{(1-x)^{k+1}(1+x)^{\lfloor(k+1)/2\rfloor}}
$$
can be used to give an expression for  $|\Omega'_{n,k}|$.

\subsection{Lucasnomials}
Sagan and Savage \cite{SaganSavage}
recently introduced analogues
of binomial coefficients, dubbed {\it lucasnomials}, defined
as follows:  for $0 \leq k \leq n$, 
$$
\left\{ \begin{matrix} n \\ k \end{matrix} \right\}:=
\frac{ \left\{ n \right\}! }
      {  \left\{ k \right\}! \left\{ n-k \right\}! }
$$
where $\left\{ n \right\}!:=
\left\{ 1 \right\} \left\{ 2 \right\} \cdots \left\{ n \right\}$,
and 
$\left\{ n \right\}$ is defined as polynomials in variables $s,t$ 
recursively, via
$
\left\{ 0 \right\}:=0 ,
\left\{ 1 \right\}:=1 ,
$
and
$$
\left\{ n \right\}
 =s \left\{ n - 1 \right\} + t \left\{ n - 2 \right\}
$$
It is not hard to see that after substituting
$$
s=q+1, \quad t=-q
$$ 
the lucasnomial is the $q$-binomial, and therefore after substituting 
\begin{equation}
\label{lucasnomial-negative-q-substitution}
s=-q+1=-(q-1), \quad t=q
\end{equation} 
the lucasnomial is (up to sign) the primed $q$-binomial
from \eqref{negative-q-binomial-definition}.
Since their main result \cite[Theorem 3.1]{SaganSavage} expands the general 
lucasnomial as a sum of monomials $s^a t^b$, 
one might wonder how their expansion
compares (after substituting as in \eqref{lucasnomial-negative-q-substitution})
with Theorem~\ref{mainth}.  It turns out that their
expansion has more terms $s^a t^b=(-1)^a (q-1)^a q^b$
than Theorem~\ref{mainth}, and not all terms in their expansion
have the same sign $(-1)^a$.

\section{Acknowledgements}
The authors thank Paul Garrett, John Shareshian, Eric Sommers, and
Ryan Vinroot for helpful comments.


\end{document}